\documentclass[11pt]{amsart}

\newtheorem{theorem}{Theorem}[section]
\newtheorem{proposition}[theorem]{Proposition}
\newtheorem{lemma}[theorem]{Lemma}
\newtheorem{cor}[theorem]{Corollary}
\newtheorem{definition}[theorem]{Definition}

\newtheorem{conjecture}[theorem]{Conjecture}

\theoremstyle{plain}
\newtheorem*{I1}{Theorem \ref{periodic}}
\newtheorem*{I2}{Theorem \ref{per height}}
\numberwithin{equation}{theorem}

\theoremstyle{remark}

\DeclareMathOperator{\Gal}{Gal}

\DeclareMathOperator{\Supp}{Supp}

\DeclareMathOperator{\Res}{Res}

\newcommand{\bP}{{\mathbb P}}
\newcommand{\bZ}{{\mathbb Z}}

\newcommand{\bC}{{\mathbb C}}

\newcommand{\bQ}{{\mathbb Q}}
\newcommand{\lra}{\longrightarrow}

\newcommand{\fO}{\frak{o}}
\newcommand{\cM}{\mathcal{M}}
\newcommand{\cU}{\mathcal{U}}

\newcommand{\bQb}{\bar{\bQ}}
\newcommand{\bK}{\overline{K}}

\title{Equidistribution and generalized Mahler measures}
\author{L.~Szpiro and T.~J.~Tucker}
\thanks{The first author was partially supported by NSF
  Grant 0071921.  The second author was partially supported by
  NSF Grant 0101636}

\address{
Lucien Szpiro \\
Ph.D. Program in Mathematics\\
Graduate Center of CUNY\\
365 Fifth Avenue\\
New York, NY 10016-4309
}

\email{lszpiro@gc.cuny.edu}

\address{
Thomas Tucker\\
Department of Mathematics\\
Hylan Building\\
University of Rochester\\
Rochester, NY 14627
}

\email{ttucker@math.rochester.edu}

\begin{document}

\begin{abstract}
  If $K$ is a number field and $\varphi: \bP^1_K \lra \bP^1_K$ is a
  rational map of degree $d > 1$, then at each place $v$ of $K$, one
  can associate to $\varphi$ a generalized Mahler measure for
  polynomials $F \in K[t]$.  These Mahler measures give rise to a
  formula for the canonical height $h_\varphi(\beta)$ of an element
  $\beta \in \bK$; this formula generalizes Mahler's formula for the
  usual Weil height $h(\beta)$. In this paper, we use diophantine
  approximation to show that the generalized Mahler measure of a
  polynomial $F$ at a place $v$ can be computed by averaging $\log
  |F|_v$ over the periodic points of $\varphi$.
\end{abstract}

\maketitle

\begin{center}
  {\it This paper is dedicated to the memory of Serge Lang, who
    taught the world number theory for more than fifty years,
    through his research, lectures, and books.}\\
\end{center}
\vspace{.5cm}


The usual Weil height of a rational number $x/y$, where $x$ and $y$
are integers without a common prime factor, is defined as $$h(x/y) =
\log \max (|x|, |y|).$$
More generally, one can define the usual Weil
height $h(\beta)$ of an algebraic number $\beta$ in a number field $K$
by summing $\log \max (|\beta|_v, |1|)$ over all of the absolute
values $v$ of $K$.  Mahler (\cite{mahler}) has proven that if $F$ is a
nonzero irreducible polynomial in $\bZ[t]$ with coprime coefficients
such that $F(\beta) = 0$, then
\begin{equation}\label{1}
 \deg(F)h(\beta) =
\int_{0}^{1} \log|F(e^{2\pi i \theta})|d\theta.
\end{equation}
The quantity $\int_{0}^{1}
\log|F(e^{2\pi i \theta})|d\theta$ is often referred to as the {\it
  Mahler measure} of $F$.

It is easy to see that $h(\beta^2) = 2 h(\beta)$ for any algebraic
number $\beta$.  Similarly, it is easy to check that for any
continuous function $g$ on the unit circle, we have
$$
\int_{0}^{1} g((e^{2\pi i \theta})^2) d\theta = \int_{0}^{1}
g(e^{2\pi i \theta}) d \theta.$$
Furthermore, the unit circle is the
Julia set of $\varphi$.  Thus, Mahler's formula says that one obtains
the height of an algebraic number by integrating its minimal
polynomial against the unique measure $\mu$ such that $\varphi^* \mu =
\mu$ and $\mu$ is supported on the Julia set of $\varphi$.

Now, let $\varphi: \bP^1_{\bC} \lra \bP^1_{\bC}$ be any nonconstant
rational map.  Brolin (\cite{Brolin}) and Lyubich (\cite{Lyubich})
have constructed a totally $\varphi$-invariant probability measure $
\mu_{\varphi}$ (that is, we have $\varphi^* \mu$ and $\varphi_* \mu$)
with support on the Julia set of $\varphi$; Freire, Lopes, and Ma\~ne
(\cite{mane}) have demonstrated that this measure is the {\it unique}
totally $\varphi$-invariant probability measure with support on the
Julia set of $\varphi$.  When $\varphi$ is defined over a number field
$K$, Call and Silverman (\cite{CS}) have constructed a height function
$h_\varphi$ with the properties that: (1) $h_\varphi(\varphi(x)) =
(\deg \varphi) h_\varphi(x)$ and (2) there is a constant $C_\varphi$
such that $|h(x) - h_\varphi (x)| < C_\varphi$ for all $x \in
\bP^1(\bK)$.  In \cite{STP}, it is shown that Mahler's formula
\eqref{1} generalizes to the adelic formula
 \begin{equation}\label{STP-mahler}
      (\deg F)  h_\varphi(x)
      = \sum_{\text{places $v$ of $K$}}\int_ {\mathbb{P}^{1}(\bC_v)}
      \log|F|_v \, d\mu_{\varphi,v},
\end{equation}
where $\beta$ is an algebraic point, $F$ is a nonzero irreducible
polynomial in $\bQ[t]$ such that $F(\beta) = 0$, the measure $
\mu_{\varphi,v}$ at an archimedean place is the totally
$\varphi$-invariant probability measure constructed by Brolin and
Lyubich, and the integral $\int_ {\mathbb{P}^{1}(\bC_v)} \log|F|_v \,
d\mu_{\varphi,v}$ at a finite place $v$ is defined so that its value
is the $v$-adic analog of the value at an archimedean place (note that
as defined in \cite{STP}, these are not integrals per se).  Favre and
Rivera-Letelier have also given a proof of \ref{STP-mahler}, using
actual integrals on Berkovich spaces; Pi{\~n}eiro (\cite{Pineiro}) and
Chambert-Loir and Thuillier (\cite{CL-T, Th}) have recently proven
higher-dimensional generalizations of \ref{STP-mahler}.

Lyubich \cite{Lyubich} has also proven that for any continuous
function $g$ and any archimedean place $v$, the integrals
$\int_{\bP^1(\bC_v)} g \, d \mu_{\varphi, v}$ can be computed by
averaging $g$ on the periodic points of $\varphi$; that is to say,
\begin{equation}\label{basic}
\lim_{k \to \infty} \frac{1}{(\deg \varphi)^k} \sum_{\varphi^k(w) = w} g(w) =
\int_{\bP^1(\bC_v)} g \, d \mu_{\varphi, v}.  
\end{equation}
Autissier (\cite{autissier}), Bilu (\cite{Bilu}), Szpiro, Ullmo, and
Zhang (\cite{suz}), and others have obtained generalizations and
variations of this result.  The most recent generalization, proven
independently by Baker and Rumely (\cite{BR}), Chambert-Loir
(\cite{CL}), and Favre and Rivera-Letelier (\cite{FR1} and \cite{FR2})
states that \eqref{basic} continues to hold when the periodic points
$w$ such that $\varphi^k(w) = w$ are replaced by the conjugates of any
infinite nonrepeating sequence of algebraic points with height tending
to 0 and when the measure $\mu_{\varphi,v}$ is the unique totally
$\varphi$-invariant measure without point masses on the $v$-adic
Berkovich space (see \cite{berkovich}) for a finite place $v$.

The function $\log |F|$, for $F$ a nonconstant polynomial, is not
continuous in general, of course.  Thus, the equidistribution results
cited above do not allow us to compute Mahler measures by averaging
$\log |F|_v$ over points of small height.  One can, however, show that
for any $\beta \in \bQb$, we have
\begin{equation}\label{circle}
[\bQ(\beta): \bQ] h(\beta) = \lim_{k \to \infty} \frac{1}{d^k}
\sum_{\xi^n = 1} \log |F(\xi)| = \int_{0}^{1} \log|F(e^{2\pi i \theta})|d\theta,
\end{equation}
where $F$ is a nonzero irreducible polynomial in $\bZ[t]$ with coprime
coefficients such that $F(\beta) = 0$ (see \cite[Chapter
1]{algebraicdynamics}, \cite{schinzel}).  Everest, Ward, and N\'{i}
Fhlath\'{u}in have proved similar results for maps that come from
multiplication on an elliptic curve (\cite[Chapter
6]{algebraicdynamics}, \cite{elliptic}).  The proofs of these results
make use of the theory of linear forms in logarithms (\cite{baker},
\cite{david}), which is used to show that the periodic points of the
maps in question have strong diophantine properties.  It is not clear
how to apply the theory of linear forms in logarithms in the case of
more general rational maps.  In this paper, we use Roth's Theorem
(\cite{Roth}) from diophantine approximation in place of the theory of
linear forms in logarithms.  This allows us to work in greater
generality.

\subsection{Statements of the main theorems}\label{statement}

The main results of this paper extend \eqref{circle} to a formula that
holds for all rational maps.  Let $K$ be a number field or a function
field of characteristic zero, let $v$ be a place of $K$, and let
$\varphi: \bP^1_K \lra \bP^1_K$ be a nonconstant rational map of
degree $d>1$.  We prove the following equidistribution result for the
periodic points of $\varphi$.

\begin{I1}
  For any nonzero polynomial $F$ with coefficients in $\bK$, we have
\begin{equation*}
\int_{\mathbb{P}^{1}(\bC_v)}
\log|F|_v \, d\mu_{\varphi,v} = 
  \lim_{{k} \to \infty} \frac{1}{d^{k}} \sum_{\substack{\varphi^{k}([w:1]) =
    [w:1] \\ F(w) \not= 0}} \log | F(w) |_v. 
\end{equation*}  
\end{I1}

This allows us to show that for any point $\beta \in \bK$, the
canonical height $h_\varphi(\beta)$ can be computed by taking the
average of the log of the absolute value of a minimal polynomial for
$\beta$ over the periodic points of $\varphi$.

\begin{I2}
  For any $\beta \in \bK$ and any nonzero
  irreducible $F \in K[t]$ such that $F(\beta) = 0$, we have
\begin{equation*}
\begin{split}
  (\deg K) & (\deg F) (h_\varphi(\beta)   - h_\varphi(\infty)) \\ 
  & =\sum_{\text{places $v$ of $K$}} \lim_{{k} \to \infty}
  \frac{1}{d^{k}} \sum_{\substack{\varphi^{k}([w:1]) = [w:1] \\ F(w)
      \not= 0}} \log | F(w) |_v.
\end{split}
\end{equation*}
\end{I2}

In both the theorems, the $w$ are counted with multiplicity.  We
explain what multiplicity means in this context in Section \ref{notation}.

We are also able to prove that $\int_{\mathbb{P}^{1}(\bC_v)} \log|F|_v
\, d\mu_{\varphi,v}$ is the limit as $n$ goes to infinity of the
average of $\log |F|_v$ on the points $w$ for which $\varphi^n(w) =
\alpha$, where $\alpha$ is an algebraic point that is not an
exceptional point for $\varphi$.  We state this in Theorem
\ref{backwards}.  This enables us to prove Theorem \ref{back height},
which is the analog of Theorem \ref{per height} for the points $w$
such that $\varphi^n(w) = \alpha$.

\subsection{Outline of the paper}\label{outline} This paper is organized as follows:
 
\ref{notation} - Notation and terminology.

\ref{brolin} - Brolin-Lyubich integrals and local heights.

\ref{dio} - Preliminaries from diophantine approximation.

\ref{main} - Main results: \ref{using} - Using Roth's Theorem;
\ref{pre} - Preperiodic points; \ref{proofs} - Proofs of the main theorems.

\ref{counterexample} - A counterexample.

\ref{applications} - Applications: \ref{lyap} - Lyapunov exponents;
\ref{symmetry pap} - Symmetry of canonical heights; \ref{small} -
Computing with points of small height.
\\
\\
\noindent The strategy of the proof of the main theorems is fairly simple.  By
additivity, it suffices to prove our results for polynomials of the
form $F(t) = t - \beta$ for $\beta \in \bK$.  After Section
\ref{brolin}, we are reduced to showing that
\begin{equation}\label{idea}
\begin{split}
\lim_{{k} \to \infty}
\frac{1}{d^{k}} \sum_{\substack{\varphi^{k}([w:1]) = [w:1] \\ w
    \not= \beta}} \log | w - \beta |_v & = \lim_{k \to \infty} \frac{\log \max
  (|P_{k}(\beta,1)|_v, |Q_{k}(\beta,1)|_v) }{d^{k}}\\
 & - \lim_{k \to \infty} \frac{\log \max
  (|P_{k}(1,0)|_v, |Q_{k}(1,0)|_v) }{d^{k}},
\end{split}
\end{equation}
where $\varphi^k$ is written as
$$\varphi^k([T_0:T_1]) = [P_k(T_0,T_1):Q_k(T_0,T_1)]$$
for coprime
homogeneous polynomials $P_k$ and $Q_k$ in the $K[T_0,T_1]$.  The
points $w$ for which $\varphi^k(w) = w$ are just the solutions to the
equation $P_k(w,1) - w Q_k(w,1) = 0$.  Thus, we get the left-hand side
of \eqref{idea} by taking the limit of $\log |P_k(\beta,1) - \beta
Q_k(\beta,1)|_v/d^k$ as $k$ goes to $\infty$.  For each $k$, we
rewrite this as
$$
\frac{\log |Q_k(\beta,1)|_v}{d^k} + \frac{\log
  |\frac{P_k(\beta,1)}{Q_k(\beta,1)} - \beta|_v}{d^k}$$
and use
diophantine approximation to show that the second term in the equation
above usually goes to 0 as $k \to \infty$; our theorems then follow
after a bit of calculation.  The diophantine approximation result we
use is Roth's Theorem, which we state in Section \ref{dio} as Theorem
\ref{Roths}.  We use Roth's Theorem to derive Lemma \ref{from Roth},
which is the key lemma in our proofs of the main theorems.  The idea
for the proof of Lemma \ref{from Roth} comes from Siegel's famous
paper \cite{Siegel}.  We should note that after writing this paper we
discovered that Silverman (\cite{SilSiegel}) has used methods very
similar to those found here at the beginning of Section \ref{main}; we
require a slight modification of his results along these lines,
however, so we present the necessary argument here in full.

Propositions \ref{bad back} and \ref{bad per} deal
with the additional complications that may arise when the $\beta$ in
\eqref{idea} is preperiodic.  These complications are overcome with
somewhat lengthy -- but essentially basic -- calculations that are
very similar to some of the computations carried out by Morton and
Silverman in \cite{MorSil2}.

In Section \ref{counterexample}, we construct a simple counterexample
that shows that Theorem \ref{periodic} will not hold in general when
the polynomial $F$ does not have algebraic coefficients (it is likely
that the theorem will also fail if the point $\alpha$ is not
algebraic).  We construct a transcendental number $\beta$ such that
the limit $\lim_{k \to \infty} \frac{1}{2^k} \sum_{\xi^{2^k} = 1} \log
|\xi - \beta|$ does not exist.  This means that there is no way to
prove the main results of this paper without using some special
properties of
algebraic numbers. \\
\\
\noindent {\it Acknowledgments.}  We would like to thank M.~Baker,
A.~Chambert-Loir, L.~DeMarco, R.~Rumely, and S.~Zhang for many helpful
conversations.  In particular, we thank M.~Baker, L.~DeMarco, and R.~Rumely for
suggesting some of the applications mentioned in Section \ref{applications}.
\ref{applications}.


    

\section{Notation and terminology}\label{notation}

We fix the following notation:
\begin{itemize}
\item $K$ is a number field or a function field of characteristic 0
  (by function field we mean a finite algebraic extension of a field
  of the form $K_{\text{cons}}(T)$ where $K_{\text{cons}}$ is
  algebraically closed in $K$);
\item $v$ is a place of $K$;
\item $K_v$ is the completion of $K$ at $v$;
\item $\bC_v$ is the completion of an algebraic closure of $K_v$ at $v$;
\item $\bK$ is the algebraic closure of $K$ in $\bC_v$ (note that this
  means that $v$ extends to all of $\bK$);
\item $n_v = [K_v: \bQ_v]$ if $K$ is a number field;
\item $n_v = 1$ if $K$ is a function field;
\item $\deg K = [K:\bQ]$ if $K$ is a number field; 
\item $\deg K = 1$ if $K$ is a function field.  
\end{itemize}

We  let $| \cdot |_v$ be an absolute value on $\bC_v$
corresponding to $v$.  When $K$ is a function field and $\pi_v$
generates the maximal prime $\cM_v$ in the local ring $\fO_v$
corresponding to $v$, we specify that
$$ | \pi_v |_v = e^{-[(\fO_v/\cM_v):K_{\text{cons}}]},$$
where $K_{\text{cons}}$ is the field of constants in $K$.  
When $K$ is a number field and $v$ is nonarchimedean, we normalize $|
\cdot |_v$ so that
$$|p|_v = p^{-n_v}$$
when $v$ lies over $p$.  When $K$ is a number
field and $v$ is archimedean we normalize so that $| \cdot |_v = |
\cdot |^{n_v}$ on $\bQ$, where $| \cdot |$ is the usual archimedean
absolute value on $\bQ$.

Throughout this paper, we will work with a nonconstant morphism
$\varphi:\bP_K^1 \lra \bP_K^1$ of degree $d >1$.  We choose
homogeneous polynomials $P,Q \in K[T_0,T_1]$ of degree $d$ without a
common factor along with a coordinate system $[s:t]$ for $\bP^1_{\bK}$ such
that
$$
\varphi([T_0:T_1]) = [P(T_0,T_1):Q(T_0,T_1)],$$
where $P$ and $Q$
have no common zero in $\bP^1(\bK)$.  We let $P_1 = P$ and $Q_1 = Q$,
and for $k \geq 2$ we define $P_k$ and $Q_k$ recursively by
$$P_k(T_0,T_1) = P_{k-1}(P(T_0,T_1), Q(T_0,T_1))$$
and 
$$Q_k(T_0,T_1) = Q_{k-1}(P(T_0,T_1),Q(T_0,T_1)).$$

Having chosen coordinates, we can define the usual Weil height as
$$
h([a:b]) = \frac{1}{\deg K}\sum_{\text{places } v \text{ of } K}
\log \max(|a|_v, |b|_v)$$
when $a,b \in K$.  When $a$ and $b$ lie in
an extension $L$ of $K$, this definition extends to
\begin{equation}\label{Weil}
h([a:b]) = \frac{1}{[L:K](\deg K)}\sum_{\text{places } w \text{ of } L}
[L_w:K_v] \log \max(|a|_w, |b|_w),
\end{equation}
where $L_w$ is the completion of $L$ at $w$ and the absolute value $|
\cdot |_w$ restricts to some $| \cdot |_v$ on $K$.

As in \cite{CS}, we define the canonical height
$h_\varphi$ 
as 
\begin{equation}\label{can}
h_\varphi([a:b]) = \lim_{k \to \infty} \frac{h(\varphi^k([a:b]))}{d^k}.
\end{equation} 

We say that $\alpha \in \bP^1(\overline{K})$ is a {\bf periodic} point
for $\varphi$ if there exists a positive integer $n$ such that
$\varphi^{n}(\alpha) = \alpha$.  If $\alpha$ is periodic, we define
the {\bf period} of $\alpha$ to be the smallest positive integer
$\ell$ such that $\varphi^{\ell}(\alpha) = \alpha$.  We say that
$\alpha$ is {\bf preperiodic} if there exists a positive integer $n$
such that $\varphi^n(\alpha)$ is periodic.

We will use a small amount of the theory of dynamics on the projective
plane; for a more thorough account of the subject, we refer the reader
to Milnor's (\cite{Milnor}) and Beardon's (\cite{Beardon}) books on
the subject.  We say that $\alpha \in \bP^1(\overline{K})$ is an {\bf
  exceptional} point for $\varphi$ if $\varphi^2(\alpha) = \alpha$ and
$\varphi^2$ is totally ramified at $\alpha$.  This is equivalent to
saying that the set $\bigcup_{k=1}^\infty (\varphi^k)^{-1} (\alpha)$
is finite (see \cite[Chapter 4.1]{Beardon}).  If $\alpha$ is
exceptional, then at each place $v$, there is a maximal $v$-adically
open set $\cU$ containing $\alpha$ such that the sequence
$(\varphi^{\ell k}(\beta))_k$ converges to $\alpha$ for each $\beta
\in \cU$, where $\ell$ is the period of $\alpha$ (which is either 1 or
2).  We call $\cU$ the {\bf attracting basin} of $\alpha$ (see
\cite[Chapter 6.3]{Beardon}, which uses the terminology ``local
basin'').

We always count points with multiplicities in this paper.  The
multiplicity of a point $[z:1]$ in the multi-set $\{ w \mid
\varphi^k(w) = w \}$ is the highest power of $t-z$ that divides the
polynomial $ P_k(t,1) - t Q_k(t,1)$.  The multiplicity of a point
$[z:1]$ in the multi-set $\{ w \mid \varphi^k(w) = [s:u] \}$ is the
highest power of $t-z$ that divides the polynomial $ u P_k(t,1) - s
Q_k(t,1)$ (here $s$, $u$, and $z$ are taken to be elements of $\bK$,
while $t$ is taken to be a variable).

We note that everything done in this paper depends upon our choice of
coordinates.  In particular, our integrals are closely related to the
canonical local heights (see \cite{CG}) for the point $[1:0]$ at
infinitely, so our choice of the point at infinity affects all of our
integrals.  To emphasize the fact that we treat $[1:0]$ as the point
at infinity, we denote it as $\infty$ where appropriate.

\section{Brolin-Lyubich integrals and local heights}\label{brolin}

We will work with the limits 
\begin{equation}\label{first lim}
\lim_{k \to \infty} \frac{\log \max (|P_{k}(a,b)|_v,
  |Q_{k}(a,b)|_v) }{d^{k}}
\end{equation}
for $(a, b) \in \bC_v \setminus \{ (0, 0) \}$.  For a proof that these
limits exist, see \cite{STP}, \cite{BR}, or \cite{CG} (the proof is
essentially an exercise in using telescoping sums and geometric
series).  Note that Call and Goldstine (\cite[Theorem 3.1]{CG}) have
shown that
$$
\hat{h}_{\varphi,v}([\beta:1]) = \lim_{k \to \infty} \frac{\log
  \max (|P_{k}(\beta,1)|_v, |Q_{k}(\beta,1)|_v) }{d^{k}}$$
is the
unique Weil function for $[1:0]$ at $v$ (see \cite[Chapter
10]{dio-geo} for a definition of Weil functions) that satisfies
$$
\hat{h}_{\varphi,v}(\varphi([a:b])) = d \hat{h}_{\varphi,v}([a:b])
+ \log \left| Q\left( \frac{a}{b},1 \right) \right|_v,$$
for any
$[a:b] \not= [1:0]$ (see \cite[Theorem 2.1]{CG}).  The function
$\hat{h}_{\varphi,v}(\cdot)$ is called a canonical local height for
$\varphi$.  We also note that these local heights can also be constructed
by taking a quantity obtained from the ``Fubini-Study'' metric and
passing to the limit; specifically, the limit in \eqref{first lim} is
also equal to
\begin{equation}\label{second lim}
\lim_{k \to \infty} \frac{\log \sqrt{|P_{k}(a,b)|_v^2 +
  |Q_{k}(a,b)|_v)^2} }{d^{k}}.
\end{equation}
The equality follows from the uniqueness of the Call-Goldstine local height
or from the arguments in \cite[Section 2]{zhangadelic}.  Note that Baker
and Rumely (\cite{BR}) use \eqref{second lim} to form local heights.

As noted in the introduction, Brolin \cite{Brolin} and Lyubich
\cite{Lyubich} have constructed a totally $\varphi$-invariant measure $
\mu_{\varphi,v}$ with support on the Julia set of $\varphi$, when $v$
is an infinite place (see also \cite{mane}).  More recently, Baker and
Rumely (\cite{BR}), Chambert-Loir (\cite{CL}), and Favre and
Rivera-Letelier (\cite{FR1} and \cite{FR2}) have constructed a
$\varphi$-invariant measure $ \mu_{\varphi,v}$ on the Berkovich space
associated to $\bP^1(\bC_v)$; this measure is unique among
$\varphi$-invariant measures without point masses the Berkovich space
associated to $\bP^1(\bC_v)$.

\begin{proposition}\label{more gen}
  Let $v$ be a place of a number field $K$ and let $F(t) = t
  - \beta$ for $\beta \in \bC_v$.  Then
\begin{equation}\label{needed}
\begin{split}
  \int_{\mathbb{P}^{1}(\bC_v)} \log|F|_v d\mu_{\varphi,v} & = \lim_{k
    \to \infty} \frac{\log \max (|P_{k}(\beta,1)|_v,
    |Q_{k}(\beta,1)|_v) }{d^{k}} \\
  & - \lim_{k
    \to \infty} \frac{\log \max (|P_{k}(1,0)|_v, |Q_{k}(1,0)|_v) }{d^{k}}.
\end{split}
\end{equation}
\end{proposition}
\begin{proof}
  We will prove this following the methods of Baker and Rumely, who
  show that the measures $\mu_{\varphi,v}$ are Laplacians of local
  height functions.  The proposition could also be proved using the
  work of Favre and Rivera-Letelier (\cite{FR2}) or Chambert-Loir and
  Thuillier (\cite{CL-T, Th}), who proved more general Mahler formulas (but
  do not formulate them in terms of limits such as \eqref{first lim}).
  In \cite{BR}, Baker and Rumely show that for $w \in \bC_v$, the
  function $H_w$ defined by
  $$
  H_w([a:b]) = - \log (wb - a) + \lim_{k \to \infty} \frac{\log
    \sqrt{|P_{k}(a,b)|_v^2 + |Q_{k}(a,b)|_v)^2} }{d^{k}}$$
  is
  subharmonic on $\bP^1(\bC_v) \setminus \{ [w:1] \}$ (see \cite{BR,
  BRbook}).  Furthermore, they show that if $\Delta$ is the
distributional Laplacian (i.e. $- d d^{c}$ considered in the
distributional sense, which can be extended to the setting of
Berkovich spaces as described in \cite{BRbook}), then
\begin{equation}\label{keyBR}
\frac{1}{n_v p(v)} \Delta H_w = - \mu_{\varphi,v} + \delta_w
\end{equation}
where $\delta_w$ is the usual Dirac point mass at $w$ and $p(v)$ is
the log of the characteristic of the residue field of $v$ when $v$
nonarchimedean, and is simply 1 when $v$ is archimedean.  Similarly,
we have
$$
\frac{1}{n_v p(v)} \Delta \log |t - \beta|_v = \delta_{[1:0]} -
\delta_\beta$$
(see \cite[Section 5.1]{FR2} or the same reasoning that
gives \eqref{keyBR}).  Now, since $\log |t - \beta|_v$ and $H_w$ are
both subharmonic on $\bP^1(\bC_v) \setminus \{ [1:0], [w:1], [\beta:1] \}$ we
have 
\begin{equation}\label{almost}
\begin{split}
\int_{\bP^1(\bC_v)} \log |t - \beta|_v d\mu_{\varphi,v} & =
\left( \int_{\bP^1(\bC_v)} \log |t - \beta|_v \left( - \frac{1}{n_v p(v)}
  \Delta H_w \right)  \right) + \log |w - \beta|_v \\
& = \left( \int_{\bP^1(\bC_v)} H_w \left ( - \frac{1}{n_v p(v)} \Delta
  \log |t - \beta|_v \right)  \right) + \log |w - \beta|_v \\
& = H_w([\beta:1]) - H_w([1:0]) + \log |w - \beta|_v
\end{split}
\end{equation}
Since \eqref{first lim} and \eqref{second lim} are equal by the
discussion above, \eqref{almost} becomes
\begin{equation*}
\begin{split}
- \log |w - \beta|_v & + \lim_{k
    \to \infty} \frac{\log \max (|P_{k}(\beta,1)|_v,
    |Q_{k}(\beta,1)|_v) }{d^{k}} + \log |1| \\
& \quad - \lim_{k
    \to \infty} \frac{\log \max (|P_{k}(1,0)|_v,
    |Q_{k}(1,0)|_v) }{d^{k}} + \log |w - \beta|_v \\
=  \lim_{k
    \to \infty} &  \frac{\log \max (|P_{k}(\beta,1)|_v, 
    |Q_{k}(\beta,1)|_v) }{d^{k}} \\
& - \lim_{k
    \to \infty} \frac{\log \max (|P_{k}(1,0)|_v,
    |Q_{k}(1,0)|_v) }{d^{k}},  
\end{split}
\end{equation*}
as desired.
\end{proof}
Note that although our integrals are defined for points in $\bC_v$,
the results we prove in Section \ref{main} only apply to points in
$\bK$.  Note as well that we make no use of the fact that our limits
correspond to actual integrals, either in the proofs of our main
theorems or in the applications in Section \ref{applications}.

When $K$ is a function field, it should also be possible to construct
suitable integrals at the places of $K$.  Since this has not yet been
done, however, we will have to make do with a definition rather than a
proof.

\begin{definition}
  Let $v$ be a place of a function field $K$ and let $F(t) = t
  - \beta$ for $\beta \in \bC_v$.  Then
\begin{equation*}
\begin{split}
  \int_{\mathbb{P}^{1}(\bC_v)} \log|F|_v d\mu_{\varphi,v} & = \lim_{k
    \to \infty} \frac{\log \max (|P_{k}(\beta,1)|_v,
    |Q_{k}(\beta,1)|_v) }{d^{k}} \\
  & - \lim_{k
    \to \infty} \frac{\log \max (|P_{k}(1,0)|_v, |Q_{k}(1,0)|_v) }{d^{k}}.
\end{split}
\end{equation*}
\end{definition}

\section{Preliminaries from diophantine approximation}\label{dio}

The following well-known theorem of Roth (\cite{Roth}) is the
principal tool from diophantine approximation that is used in this
paper.
\begin{theorem}\label{Roths} (Roth).
If $\alpha \in \bC$ is algebraic over $\bQ$,
then for any $\epsilon > 0$, there is a constant $C$ such that
\begin{equation*}
\left| \alpha - \frac{a}{b} \right| > \frac{C}{|b|^{2 + \epsilon}},
\end{equation*}
for all $a/b \in \bQ$ such that $a/b \not= \alpha$.  
\end{theorem}


We will need to work in slightly greater generality.  In the
terminology of the previous section, Roth's admits the following
generalization (see \cite[Theorem 7.1.1] {dio-geo}), which holds when
$K$ is number field or a function fields of characteristic 0.

\begin{theorem}\label{roth2}
  Let $\alpha_1, \dots, \alpha_n$ be elements of ${\overline K}$ and
  let $L \subset \bK$ be a finite extension of $K$.  Then, for any
  $\epsilon > 0$ and any places $v$ of $K$ and $w$ of $L$ such that $w
  | v$, we have
  $$
  \frac{1}{[L:K] (\deg K)} \sum_{i=1}^n \max(0, - \log| \alpha_i - \beta
  |^{[L_w:K_v] n_v}_v) \leq (2 + \epsilon) h(\beta) + O(1),$$
  for all
  $\beta \in L$ not in the set $\{ \alpha_1, \dots, \alpha_n\}$.
\end{theorem}

Let $[a:1]$ be a point in $\bP^1(\bK)$.  Then for any $[b:1] \not=
[a:1]$ in $\bP^1(\bC_v)$, we let
$$
\lambda_{[a:1],v}([b:1]) = \max(- \log |b - a|_v, 0).$$
We extend this definition to the point at $[1:0]$ by letting
$$ \lambda_{[a:1],v}([1:0]) = 0.$$
and
\begin{equation}\label{in-def}
\lambda_{[1:0],v}([b:1]) = \max(0, \log |b|_v).
\end{equation}

We will work with divisors on $\bP_{\overline K}^1$ rather than
elements of ${\overline K}$.  Let $D = \sum_{i=1}^n m_i \alpha_i$,
where $\alpha_i \in \bP^1(\bK)$ and $m_i \in \bZ$.  We let
$$\lambda_{D,v}(\beta) = \sum m_i \lambda_{\alpha_i, v}(\beta)$$
for
points $\beta \in \bP^1(\bC_v)$ that are not in $\Supp D$.  Then
$\lambda_{D,v}$ is a {\bf Weil function} for $D$ at $v$ as defined in
\cite[Chapter 10]{dio-geo}.  It is easy to check that for any divisor
$D$ and any rational map $\varphi$ on
$\bP^1$, we have 
\begin{equation}\label{functorial}
\lambda_{D,v} (\varphi(\beta)) = \lambda_{\varphi^*D, v}(\beta) + O(1),
\end{equation}
for all $\beta \in \bP^1(\bK)$ away from the support of $D$
and $\varphi^* D$.  This is a general functorial property of Weil
functions, as explained in \cite[Chapter 10]{dio-geo}.

For a divisor $D = \sum_{i=1}^n m_i \alpha_i$, where
  $\alpha_i \in \bP^1(\overline K)$, we
  define
  $$
  r(D) = \max_i(m_i).$$
  With this terminology, it follows from
  Theorem \ref{roth2} that for any $\epsilon > 0$, any finite
  extension $L$ of $K$, and any positive divisor $D$ on
  $\bP^1({\overline K})$ with $r(D) = 1$, we have
$$\frac{1}{[L:K] (\deg K)} \lambda_{D,v} (\beta) \leq (2 + \epsilon)
h(\beta) + O(1)$$
for all $\beta \in \bP^1(L)$ away from the support
of $D$.  Hence, for any positive divisor $D$ we have
\begin{equation} \label{Roth}
\frac{1}{[L:K] (\deg K)}\lambda_{D,v} (\beta)
  \leq r(D) (2 + \epsilon) h(\beta) + O(1).
\end{equation}

\section{Main results}\label{main}


We begin with a simple Lemma on how $r((\varphi^{n})^*(D)$ behaves as $n
\to \infty$ when $D$ is a divisor that does not contain an
exceptional point of $\varphi$.  We recall that in general if $D =
\sum_{i=1}^n m_i \alpha_i$ is a
divisor on $\bP^1$ and $\psi: \bP^1 \lra \bP^1$ is a nonconstant
rational map, then 
\begin{equation}\label{rami}
\psi^* D = \sum_{i=1}^n \sum_{\psi(\beta_i) = \alpha_i}
m_i e(\beta_i/\alpha_i) \beta_i
\end{equation}
where $e(\beta_i/\alpha_i)$ is the ramification
index of $\psi$ at $\beta_i$.

\begin{lemma}\label{r(D)}
Let $D$ be a divisor such that $\Supp D$ does not contain any
exceptional points of $\varphi$.  Then
$\lim_{k
  \to \infty} \frac{r((\varphi^k)^*D)}{d^k} = 0.$
\end{lemma}
\begin{proof}
  Recall that $\alpha$ is an exceptional point if and only if
  $\varphi^2(\alpha) = \alpha$ and $\varphi$ is totally ramified at
  both $\alpha$ and $\varphi(\alpha)$.  Since $\varphi$ has at most
  two totally ramified points, it follows that if $\alpha$ is not
  exceptional, then one of $\alpha$, $\varphi(\alpha)$, and
  $\varphi^2(\alpha)$ is not a totally ramified point of $\varphi$.
  Since the degree of $\varphi^3$ is $d^3$, this means that for any
  divisor $E$ such that $\Supp E$ does not contain an exceptional
  point, we have $r((\varphi^3)^* E) < d^3 r(E)$ (by \eqref{rami}), so
  $r((\varphi^3)^* E) \leq d^2(d-1) r(E)$ Now, since $\Supp D$ does
  not contain an exceptional point, $\Supp (\varphi^k)^*D$ does not
  contain an exceptional point for any $k$.  Thus, for any $k \geq 3$,
  we see that $\frac{r((\varphi^k)^*D)}{d^k}$ is less than or equal to
  $((d-1)/d)^{(k-2)/3} r(D),$ which goes to zero as $k$ goes to
  infinity.
\end{proof}

\subsection{Using Roth's Theorem}\label{using} Roth's Theorem allows us to prove
the following lemma.  The idea of the proof is that if $\varphi^{k
  +\ell}(\beta)$ approximates $D$ very closely, then
$\varphi^k(\beta)$ approximates $(\varphi^\ell)^*D$ very closely.
Since $\varphi^k(\beta)$ has height approximately equal to $1/d^\ell$
times the height of $\varphi^{k +\ell}(\beta)$, this makes
$h(\varphi^k(\beta))$ small relative to
$\lambda_{(\varphi^\ell)^*D}(\beta)$.  Repeating this for infinitely
many $\varphi^k(\beta)$ gives a contradiction to Roth's Theorem.  This
idea is due to Siegel (\cite{Siegel}); similar arguments can be found
in \cite{SilSiegel}.
 
\begin{lemma}\label{from Roth}
  Let $D$ be a positive divisor on $\bP^1$ such that $\Supp D$ does
  not contain any of the exceptional points of $\varphi$.  Let $\beta$
  be a point in $\bP^1(\bK)$ for which there is a strictly increasing
  sequence of integers $( e_i )_{i=1}^\infty$ such that
  $\varphi^{e_i}(\beta) \notin \Supp D$.  Then
\begin{equation}\label{key}
\lim_{i \to \infty} \frac{\lambda_{D,
    v}(\varphi^{e_i}(\beta))}{d^{e_i}} = 0.
\end{equation}
\end{lemma}
\begin{proof}
  
  Let $L$ be a finite extension of $K$ for which $\beta \in \bP^1(L)$.
  Choose $\delta > 0$.  By Lemma \ref{r(D)}, we may pick an integer
  $\ell$ such that $\frac{r((\varphi^\ell)^*D)}{d^\ell} < \delta/2$.
  We may then write $ \frac{r((\varphi^{\ell})^*D) (2 +
    \epsilon)}{d^\ell} = \delta$ for some $\epsilon > 0$.  For any
  $e_i$, we have $\varphi^{e_i - \ell} (\beta) \notin \Supp
  (\varphi^\ell)^*D$ since $\varphi^{e_i} (\beta) \notin \Supp D$.
  Thus, applying Roth's Theorem (as expressed in \eqref{Roth}), we find
  that for all $e_i$ we have
  $$
  \frac{1}{[L:K] (\deg K)} \lambda_{(\varphi^{\ell})^*D, v} (\varphi^{e_i
    - \ell}(\beta)) \leq r((\varphi^{\ell})^*D) (2 + \epsilon)
  h(\varphi^{e_i - \ell} (\beta)) + O(1).$$  Using
  \eqref{functorial} and the fact that 
  $h(\varphi^{e_i} (\beta)) \leq d^\ell
  h_\varphi(\varphi^{e_i - \ell}(\beta)) + O(1)$, we then obtain 
\begin{equation*}
\begin{split}
  \frac{1}{[L:K] (\deg K)} \lambda_{D,v} (\varphi^{e_i}(\beta)) & \leq
  \frac{1}{[L:K] (\deg K)} \lambda_{(\varphi^\ell)^*D,
    v} (\varphi^{e_i - \ell}(\beta)) + O(1)\\
  & \leq r((\varphi^{\ell})^*D) (2 +
  \epsilon) h(\varphi^{e_i - \ell} (\beta)) + O(1) \\
& \leq  \frac{r((\varphi^{\ell})^*D) (2 +
  \epsilon)}{d^\ell} h(\varphi^{e_i}(\beta)) +O(1)\\
  & \leq \delta h(\varphi^{e_i}(\beta)) +O(1)\\
  & \leq \delta d^{e_i} h(\beta) + O(1).
\end{split} 
\end{equation*}
Dividing through by $d^{e_i}$ gives
$$
\lim_{i \to\infty} \sup \frac{\lambda_{D,v}
  (\varphi^{e_i}(\beta))}{d^{e_i}} \leq [L:K] (\deg K) \delta
h(\beta).$$
Since $\lambda_{D,v} (\varphi^{e_i}(\beta)) \geq 0$,
letting $\delta$ go to zero gives \eqref{key}, as desired.
\end{proof}

This allows us to prove the following Proposition, which will be used
to prove Theorems \ref{backwards} and \ref{periodic}.
 
\begin{proposition}\label{to use}
  Let $\alpha = [s:u]$ be a nonexceptional point in $\bP^1(\bK)$.
  Then for any point $\beta = [a:b]$ in $\bP^1(\bK)$ and any strictly
  increasing sequence of integers $( e_i )_{i=1}^\infty$ such that
  $\varphi^{e_i}(\beta) \not= \alpha$, we have
$$
\lim_{i \to \infty} \frac{\log |uP_{e_i}(a,b) - s Q_{e_i}(a,b)|_v}{d^{e_i}} =
\lim_{i \to \infty} \frac{\log \max (|P_{e_i}(a,b)|_v,
  |Q_{e_i}(a,b)|_v) }{d^{e_i}}. $$
\end{proposition}
\begin{proof}
  Note that we know that the the limit on the right-hand side of the
  equation above exists by the discussion at the beginning of
  Section~\ref{brolin}.  
  
  If $[1:0]$ is an exceptional point of $\varphi$, let $\cU$ be its
  attracting basin; if $[1:0]$ is not exceptional, let $\cU$ simply
  equal $\{ [1:0] \}$.  We will divide $( e_i )_{i=1}^\infty$ into two
  subsequences: one consisting of the $e_i$ for which
  $\varphi^{e_i}(\beta) \notin \cU$ and one consisting of the
  remaining integers in the sequence $( e_i )_{i=1}^\infty$.  Let
  $(\ell_j)_{j=1}^\infty$ be the subsequence consisting of all integers
  $\ell_j$ in $( e_i )_{i=1}^\infty$ such that $\varphi^{\ell_j}(\beta)
  \notin \cU$ (this subsequence may be empty).  We have
  \begin{equation} \label{zero}
  \begin{split}
  \lim_{j \to \infty} \frac{ \max(\log
  |P_{\ell_j}(a,b)/Q_{\ell_j}(a,b)|_v, 0)}{d^{\ell_j}} = 0.
  \end{split}  
  \end{equation}
  If $[1:0]$ is not exceptional, this follows from Lemma \ref{from
    Roth} applied to $D = [1:0]$, along with \eqref{in-def}.  If
  $[1:0]$ is exceptional, the fact that $\varphi^{\ell_j}(\beta)
  \notin \cU$ for all $j$ implies that
  $|P_{\ell_j}(a,b)/Q_{\ell_j}(a,b)|_v$ is bounded for all $j$, so
  \eqref{zero} clearly holds.  It follows immediately from
  \eqref{zero} that
  \begin{equation}\label{Qmax}
  \lim_{j \to \infty} \frac{\log \max (|P_{\ell_j}(a,b)|_v,
    |Q_{\ell_j}(a,b)|_v) }{d^{\ell_j}} = \lim_{j \to \infty} \frac{\log
    |Q_{\ell_j}(a,b)|_v}{d^{\ell_j}}. 
\end{equation}
  Note that if $u = 0$, then
  $$uP_{\ell_j}(a,b) - sQ_{\ell_j}(a,b) = s Q_{\ell_j}(a,b),$$ 
  so we are done.
  Otherwise, by Lemma \ref{from Roth}, we have
  $$
  \lim_{j \to \infty} \frac{ \max \left( 0, - \log
    \left| \frac{P_{\ell_j}(a,b)}{Q_{\ell_j}(a,b)} -
      \frac{s}{u} \right|_v \right)} {d^{\ell_j}} = 0.$$
Combining this with \eqref{zero}, we see that
$$
\lim_{j \to \infty} \frac{\log
   \left| \frac{P_{\ell_j}(a,b)}{Q_{\ell_j}(a,b)} - \frac{s}{u} \right|_v}
{d^{\ell_j}} = 0.$$

Thus, using \eqref{Qmax}, we obtain
\begin{equation*}
\begin{split}
  \lim_{j \to \infty} & \frac{\log |uP_{\ell_j}(a,b) - s
    Q_{\ell_j}(a,b)|_v}{d^{\ell_j}} \\ & = \lim_{j \to \infty}
  \frac{\log \left( |Q_{\ell_j}(a,b)|_v |u|_v \left|
        \frac{P_{\ell_j}(a,b)}{Q_{\ell_j}(a,b)} -
        \frac{s}{u} \right|_v \right)}{d^{\ell_j}} \\
  & = \lim_{j \to \infty} \frac{\log |Q_{\ell_j}(a,b)|_v}{d^{\ell_j}}
  + \lim_{j \to \infty} \frac{ \log \left|
      \frac{P_{\ell_j}(a,b)}{Q_{\ell_j}(a,b)} -
      \frac{s}{u} \right|_v}{d^{\ell_j}} \\
  & = \lim_{j \to \infty} \frac{\log \max (|P_{\ell_j}(a,b)|_v,
    |Q_{\ell_j}(a,b)|_v) }{d^{\ell_j}},
\end{split}
\end{equation*}
as desired.  

Now, let $(m_j)_{j=1}^\infty$ be be the subsequence of
$(e_i)_{i=1}^\infty$ consisting of all integers $m_j$ in
$(e_i)_{i=1}^\infty$ such that $\varphi^{m_j}(\beta) \in \cU$ (this
subsequence may also be empty).  If $\alpha = [1:0]$, then $[1:0]$ is
not exceptional by assumption, so there are no $m_j$ and we are done.
Otherwise, we have
$$\lim_{j \to \infty} \frac{|s Q_{m_j} (a,b)|_v}{|u P_{m_j}(a,b)|_v}
= 0, $$
since $\frac{P_{m_j}(a,b)}{Q_{m_j}(a,b)}$ goes to infinity and $u
\not= 0$.  This implies that
\begin{equation*}
\begin{split}
  & \lim_{j \to \infty} \frac{\log |uP_{m_j}(a,b) - s
    Q_{m_j}(a,b)|_v}{d^{m_j}}
\\  & =
  \lim_{j \to \infty} \frac{\log |uP_{m_j}(a,b)|_v}{d^{m_j}} \\
  & = \lim_{j \to \infty} \frac{\log \max (|P_{m_j}(a,b)|_v,
    |Q_{m_j}(a,b)|_v)) } {d^{m_j}}.
\end{split}
\end{equation*}

Since every element of the sequence $( e_i )_{i=1}^\infty$ is in $(
\ell_j )_{j=1}^\infty$ or $(
m_j )_{j=1}^\infty$, this completes our proof.  
\end{proof}

\subsection{Preperiodic points}\label{pre}
Proposition \ref{to use} provides all the information we need when
$\varphi^k([a:b]) = [s:u]$ for at most finitely many $k$; this will
always be the case when $[s:u]$ is not preperiodic.  When $[s:u]$ is
preperiodic, however, there may be infinitely many $k$ such that
$\varphi^k([a:b]) = [s:u]$.  New complications arise when this is the
case; we treat these complications in Propositions \ref{bad back} and
\ref{bad per}.

Suppose that $(bT_0 - aT_1)^{w_k}$ is the highest power of $(bT_0 -
aT_1)$ that divides $u P_k(T_0,T_1) - s Q_k(T_0,T_1)$ in
$\bK[T_0,T_1]$.  We write
$$ u P_k(T_0,T_1) - s Q_k(T_0,T_1) = (bT_0 - aT_1)^{w_k}
G_k(T_0,T_1)$$
where $G_k$ is a polynomial in $\bK[T_0,T_1]$ such that
$G_k(a,b) \not= 0$.

\begin{proposition}\label{bad back}
  Let $[s:u]$ be a nonexceptional point of $\varphi$.  Then, with
  notation as above, we have
\begin{equation}\label{bad back eq}
  \lim_{k \to \infty} \frac{\log |G_k(a,b)|_v }{d^k} = \lim_{k \to
    \infty} \frac{\log \max (|P_{k}(a,b)|_v, |Q_{k}(a,b)|_v )}{d^{k}}.
\end{equation}
\end{proposition}
\begin{proof}
  By Proposition \ref{to use}, equation \eqref{bad back eq} holds if
  we restrict to the $k$ for which $\varphi^k ([a:b]) \not= \alpha$.
  If there are only finitely many $k$ such that $\varphi^k ([a:b])
  = \alpha$, we are therefore done.  Otherwise, let $j$ be the
  smallest positive integer such that $\varphi^j([\beta:1]) = \alpha$
  and let $\ell$ be the period of $\alpha$.  Then
  $\varphi^k([\beta:1]) = \alpha$ precisely when $k$ is of the form $j
  + m \ell$ for some integer $m \geq 0$.  If $\varphi^{\ell} ([s:u]) =
  [s:u]$, then $u T_0 - s T_1$ divides $u P_{\ell}(T_0,T_1) - s
  Q_{\ell}(T_0,T_1)$.
  
  Suppose that $u \not= 0$.  Then, expanding $Q_\ell$ out in the
  variables $uT_0 - sT_1$ and $T_1$, we see that since $uT_0 - sT_1$
  cannot divide $Q_\ell(T_0, T_1)$ (because if it did, then it would
  also divide $P_\ell(T_0, T_1)$ and we know that $Q_\ell$ and
  $P_\ell$ have no factors), we have
$$ Q_\ell (T_0,T_1) = g_0 T_1^{d^\ell} + (uT_0 - sT_1) W(T_0,T_1)$$
for some nonzero $g_0 \in \bK$ and some $W(T_0,T_1) \in \bK[T_0,T_1]$.
For any $m \geq 1$ we thus have 
\begin{equation*}
  Q_{m \ell} = g_0 (Q_{(m-1)\ell})^{d^\ell} 
   + (uP_{(m-1) \ell} - s Q_{(m-1) \ell})
  W(P_{(m-1) \ell}, Q_{(m-1) \ell}).
\end{equation*}
Using induction, we see then that
\begin{equation}\label{Q}
Q_{m \ell}(T_0,T_1) = g_0^{\sum_{i=0}^{m-1} d^{i \ell}}
 T_1^{d^{m \ell}} +
(uT_0 - sT_1) W_m(T_0,T_1),
\end{equation}
for some polynomial $W_m(T_0,T_1) \in \bK[T_0,T_1]$.   
Similarly, we may write 
\begin{equation}\label{r}
\begin{split}
& u P_{\ell}(T_0,T_1) - s Q_{\ell}(T_0,T_1) \\
& = (u T_0 - s T_1)^r f_r
T_1^{d-r} + (u T_0 - s T_1)^{r+1} V(T_0,T_1),
\end{split}
\end{equation}
for some nonzero $f_r \in \bK$, some integer $r >0$, and some
$V(T_0,T_1)$ in $\bK[T_0,T_1]$.  Since $[s:u]$ is not an exceptional
point of $\varphi$, we have $r < d^{\ell}$ (note that if $r$ were to
equal to $d^{\ell}$, then $\varphi$ would have to ramify totally at
$\varphi([s:u]), \dots, \varphi^\ell{[s:u]}$, which would imply that
$\ell = 2$ and that $[s:u]$ is therefore an exceptional point, as
explained in Section \ref{notation}).  Then for any $m$, we have
\begin{equation*}
\begin{split}
  u P_{m \ell} - s Q_{m \ell} & = (u P_{(m-1) \ell} -
  s Q_{(m-1)\ell})^r f_r Q_{(m-1) \ell} ^{d-r} \\
  & + (P_{(m-1)\ell} - s Q_{(m-1)\ell})^{r+1}
  V(P_{(m-1)\ell},Q_{(m-1)\ell}),
\end{split}
\end{equation*}
so, using \eqref{Q}, \eqref{r}, and induction, we obtain
\begin{equation}\label{d-r}
\begin{split}
  & u P_{m \ell}(T_0,T_1) - s Q_{m \ell}(T_0,T_1) \\
  & = (uT_0 - s T_1)^{r^m} f_r^{\sum_{i=0}^{m-1} r^i} T_1^{d^{m \ell}
    - r^{m}} g_0^{\sum_{i=0}^{m-1} (d^{i \ell} - r^i)}\\
& + (uT_0 - s T_1)^{r^m+1} Z_m(T_0,T_1),
\end{split}
\end{equation}
for $Z_m$ a polynomial in $\bK[T_0,T_1]$.  
Since $r < d^\ell$, we have
$$ \lim_{m \to \infty} \frac{\log |f_r^{\sum_{i=0}^{m-1} r^i}
g_0^{\sum_{i=0}^{m-1} (d^{i \ell} - r^i)} |_v}{d^{m \ell}}  =
\lim_{m \to \infty} \frac{\log| g_0^{\sum_{i=0}^{m-1} d^{i \ell}} |_v}{d^{m
    \ell}} = \frac{\log |g_0|_v}{d^{\ell} - 1}.$$

Now, let $\epsilon$ be the highest power of $a T_0 - bT_1$ that divides $u
P_j - s Q_j$.  Using \eqref{d-r}, we see that we have
$$
u P_{j+m \ell}(T_0,T_1) - s Q_{j + m \ell}(T_0,T_1) = (bT_0 -
aT_1)^{\epsilon r^m} G_{j + m \ell}(T_0,T_1)$$
for a polynomial $G_{j +
  m \ell} \in \overline{K}[T_0,T_1]$.  Letting $m$ go to infinity, we
see from \eqref{d-r} that
\begin{equation*}
  \lim_{m \to \infty} \frac{\log|G_{j + m \ell}(a,b)|_v} {d^{j + m \ell}} 
    = \frac{\log|g_0|_v}{d^j(d^{\ell} - 1)} + \frac{\log|Q_j(a,b)|_v}{d^j}.
\end{equation*}  

Similarly, \eqref{Q} yields
$$
\lim_{m \to \infty} \frac{\log |Q_{j + m \ell}(a,b)|_v}{d^{j + m
    \ell}} = \frac{\log| g_0|_v}{d^j(d^{\ell} - 1)} + \frac{\log |Q_j(a,b)|_v}{d^j}.$$

Moreover, since 
$ u P_{j + m \ell}(a,b) = s Q_{j + m \ell}(a,b)$ for every $m$, we
have 
$$
\lim_{m \to \infty}\frac{\log |P_{j + m \ell}(a,b)|_v}{d^{j + m
    \ell}} = \lim_{m \to \infty}\frac{\log |Q_{j + m \ell}(a,b)|_v}{d^{j + m
    \ell}}.$$
Hence 
$$
\lim_{m \to \infty} \frac{\log| G_{j + m \ell}(a,b)|_v}{d^{j + m \ell}} =
\lim_{m \to \infty} \frac{\log \max (|P_{j + m \ell}(a,b)|_v, |Q_{j +
    m \ell}(a,b)|_v ) }{d^{j + m \ell}},$$
which completes our proof in
the case $u \not= 0$.  The proof in the case $u = 0$ proceeds in
exactly the same way, using $T_0$ in place of $T_1$.

\end{proof}

We have a similar result for the polynomials $T_0 P_k - T_1 Q_k$.  We
write
$$
T_0 P_k(T_0,T_1) - T_1 Q_k(T_0,T_1) = (bT_0 - aT_1)^{n_k}
H_k(T_0,T_1)$$
where $H_k$ is a polynomial in $\bK[T_0,T_1]$ such that
$H_k (a,b) \not= 0$.  The proof of the following proposition is
similar to Morton's and Silverman's proof of \cite[Lemma
3.4]{MorSil2}, but it requires a bit more detail since it yields
information about $H_k(a,b)$ as well as $n_k$.

\begin{proposition}\label{bad per}
With notation as above, we have
\begin{equation}\label{bad per eq}
\lim_{k \to \infty} \frac{\log |H_k(a,b)|_v }{d^k} = \lim_{k \to
  \infty} \frac{\log \max (|P_{k}(a,b)|_v, |Q_{k}(a,b)|_v) }{d^{k}}.
\end{equation}
Furthermore, $n_k$ remains bounded as $k$ goes to infinity. 
\end{proposition}
\begin{proof}
If  $( e_i )_{i=1}^\infty$ is a strictly increasing sequence
  of integers such that
  $\varphi^{e_i}([a:b]) \not= [a:b]$ for each $e_i$, then 
  $$
  H_{e_i} (T_0,T_1) = T_0 P_{e_i}(T_0,T_1) - T_1 Q_{e_i}(T_0,T_1)$$
  for all $e_i$.  Hence, by Proposition \ref{to use}, we
  $$
  \lim_{i \to\infty} \frac{\log |H_{e_i}(a,b) |_v}{d^{e_i}} =
  \lim_{i \to \infty} \frac{\log \max (|P_{e_i}(a,b)|_v,
    |Q_{e_i}(a,b)|_v) }{d^{e_i}}. $$
  If $[a:b]$ is not periodic, this
  finishes the proof.  Thus, we may assume that $[a:b]$ is periodic.  The
  rest of the proof is a computation.  We divide it into three steps.\\
 \\ 
\noindent  {\bf Step I.} We begin by changing variables so that $[a:b]$ becomes
$[0:1]$.  If $b = 0$, we write $U_0 = T_1/a $ and $U_1 = - T_0$.  We
then let
$$ R(U_0,U_1) = \frac{1}{a} Q(T_0,T_1)$$
and 
$$S(U_0,U_1) = - P(T_0,T_1)$$
(this is simply the inverse of the
transformation we defined on $T_0$ and $T_1$ -- our change of
variables is obtained by conjugation by a change-of-basis matrix).  If
$b \not= 0$, we write $U_1 = \frac{1}{b} T_1$ and
$$ U_0 = b T_0 - a T_1.$$
We then let $S(U_0,U_1) =  Q(T_0,T_1)/b$ and 
$$ R(U_0,U_1) = b P(T_0,T_1)  - a Q(T_0,T_1).$$
We define $R_k$ and $S_k$ recursively by letting $R_1 = R$, $S_1 = S$,
and setting 
$$R_{k+1}(U_0,U_1) = R_k(R(U_0,U_1), S(U_0,U_1))$$
and
$$S_{k+1}(U_0,U_1) = S_k(R(U_0,U_1), S(U_0,U_1)).$$
By the construction
of our change of variables, we have
\begin{equation}\label{same}
U_1 R_k(U_0, U_1) - U_0 S_k(U_0, U_1) = T_0 P_k(T_0,T_1) - T_1
Q_k(T_0,T_1)
\end{equation}
as polynomials in $T_0$ and $T_1$.  Hence, if $U_0^{n_k}$ is the
highest power of $U_0$ that divides $U_1 R_k(U_0, U_1) - U_0 S_k(U_0,
U_1)$ and $\tau_{k}$ is the coefficient of the $U_0^{n_k}U_1^{d^k -
  n_k}$ term in $U_1 R_k(U_0, U_1) - U_0 S_k(U_0, U_1)$, then
$$ \tau_k = H_k(a,b).$$ 

Now, let $\ell$ be the smallest positive integer for which
$\varphi^{\ell}([a:b]) = [a:b]$.  Note that $|S_{m \ell} (1,0)|_v =
\frac{|Q_{m \ell} (a,b)|_v}{|b|_v}$ if $b \not= 0$ and $ |S_{m
  \ell}(1,0)|_v = |P_{m \ell}(a,b)|_v/|a|_v$ otherwise.  Since
$$ [ P_{m \ell}(a,b):  Q_{m \ell}(a,b)] = [a:b]$$
for every $m$, it follows that  
$$
\lim_{m \to \infty} \frac{ \log|S_{m \ell}(0,1)|_v }{d^{m \ell}} =
\lim_{m \to \infty} \frac{\log \max (|P_{m \ell}(a,b)|_v, |Q_{m
    \ell}(a,b)|_v) }{d^{m \ell}}.$$ 
Thus, it will suffice to show that
\begin{equation}\label{as desired}
 \lim_{m \to \infty} \frac{ \log|\tau_{m \ell}|_v}{d^{m \ell}} = \lim_{m
  \to \infty} \frac{ \log|S_{m \ell}(0,1)|_v }{d^{m \ell}}.
\end{equation}
We write
$$ R_{\ell} (U_0, U_1) = \sum_{i=1}^{d^\ell} f_i U_0^i U_1^{d^{\ell} - i} $$
(note that $U_0$ divides $R_\ell$ by our change of variables) and
$$ S_{\ell} (U_0, U_1) = \sum_{i=0}^{d^\ell} g_i U_0^i U_1^{d^{\ell} - i}.$$
Using induction, we see that 
$$
R_{m \ell}(U_0,U_1) \equiv f_1^m g_0^{(\sum_{j=0}^{m-1} d^{ j \ell}) -
  m} U_0 U_1^{d^{m \ell} - 1} \pmod{U_0^2} $$
and
$$ S_{m \ell}(U_0,U_1) \equiv  g_0^{\sum_{j=0}^{m-1} d^{ j \ell}} U_1^{d^{m
    \ell} } \pmod{U_0^2}. $$
Thus, we have
\begin{equation}\label{U}
\begin{split}
& U_1 R_{m \ell}(U_0,U_1) - U_0 S_{m \ell} (U_0, U_1) \\
& \equiv g_0^{\sum_{j=0}^{m-1} d^{ j \ell}}
((f_1/g_0)^m - 1 )  U_0 U_1^{d^{m \ell} } \pmod{U_0^2}. 
\end{split}
\end{equation}

\noindent {\bf Step II.}  We will now treat the $m$ for which $(f_1/g_0)^m \not= 1$
We have
$$
|\log |(f_1/g_0)^m - 1 |_v \leq h((f_1/g_0)^m - 1 ) \leq 2 m
[K(f_1/g_0):K] h(f_1/g_0)$$
for all $m$ such that $(f_1/g_0)^m \not=
1$ (this is a simple version of Liouville's theorem), so
$$
\lim_{\substack {m \to \infty \\ (f_1/g_0)^m \not= 1 }} \frac{\log
    |(f_1/g_0)^m - 1  |_v}{d^{m \ell}} = 0.$$
Thus, dividing \eqref{U} through by $U_0$, we obtain
\begin{equation*}
  \lim_{\substack {m \to \infty \\ (f_1/g_0)^m \not= 1 }}
  \frac{\log |\tau_{m \ell}|_v}{d^{m \ell}}  = \lim_{m \to \infty} \frac{ \log
  |g_0^{\sum_{j=0}^{m-1} d^{ j \ell}}|_v}   {d^{m \ell}} 
    = \lim_{m \to \infty} \frac{\log |S_{m \ell}(0,1)|_v}{d^{m
        \ell}},
\end{equation*}
as desired.\\
\\
\noindent {\bf Step III.} We are left with treating the $m$ for which
$(f_1/g_0)^m = 1$.  Let $\rho$ be the smallest positive integer $m$
such that $ (f_1/g_0)^m = 1$ and write $\omega = \rho \ell$.  For
$q \geq 1$ we write
$$
R_{q \omega}(U_0,U_1) = \sum_{i=1}^{d^{q \omega}} x^{[q]}_i U_0^i
U_1^{d^{q \omega} - i}$$
(the summation starts at 1 since $U_0$
divides $R_{q \omega}$) and
$$
S_{q \omega}(U_0,U_1) = \sum_{i=0}^{d^{q \omega}} y^{[q]}_i U_0^i
U_1^{d^{q \omega} - i}.$$

Since $f_1^{\rho} = g_0^\rho$ by assumption, we have $y^{[1]}_0 =
x^{[1]}_1$ by \eqref{U}.  Multiplying $R_\omega$ and $S_\omega$ through by a
constant will change all of the limits we are calculating by the same
fixed amount, so we may assume that $y^{[1]}_0 = x^{[1]}_1 = 1$.  Let
$r$ be the smallest integer greater than 0 such that $x^{[1]}_r \not=
y^{[1]}_{r-1}$ (we have $r \geq 2$ since $(f_1/g_0)^m =1$).  Then
$U_0^r$ divides $U_1 R_{\omega} - U_0 S_{\omega}$, which in turn
divides $U_1 R_{q \omega} - U_0 S_{q \omega}$ for any $q$; hence
$U_0^r$ divides $U_1 R_{q \omega} - U_0 S_{q \omega}$ for every $q$,
so $x^{[q]}_j = y^{[q]}_{j-1}$ for $j < r$.  To calculate $x^{[q]}_r -
y^{[q]}_{r-1}$, we introduce some notation: we let
$$\left( \sum_{i=0}^M t_i U_0^i U_1^{M-i} \right)_j  = t_j$$
 for any
polynomial $\sum_{i=0}^M t_i U_0^i U_1^{M-i}$.  We have
\begin{equation}\label{sum thing}
\begin{split}
  & x^{[q]}_r - y^{[q]}_{r-1} \\
& = \sum_{i=1}^r x^{[q-1]}_i \left(
    (R_\omega)^i (S_\omega)^{d^{(q-1)\omega} - i} \right)_{r} -
  \sum_{j=0}^{r-1} y^{[q-1]}_j \left( (R_\omega)^j (S_\omega)^{d^{(q-1)
        \omega} - j} \right)_{r - 1}.
\end{split}
\end{equation}
For any $i<r$, we have $x^{[1]}_i = y^{[1]}_{i-1}$, so $(U_0
R_\omega)_i = (U_1 S_\omega)_i$.  Hence, we have
$$ \left( (R_\omega)^j (S_\omega)^{d^{(q-1)
        \omega} - j} \right)_{r - 1} = \left( (R_\omega)^{j+1} (S_\omega)^{d^{(q-1)
        \omega} - j -1} \right)_r$$
for $j > 0$.  
For $j=0$, we have
\begin{equation*}
\begin{split}
  \left( S_\omega^{d^{(q-1)\omega}} \right)_{r - 1} & = \left( \big( R_\omega +
    (x^{[1]}_r - y^{[1]}_{r-1})U_0^r U_1^{d^{\omega} - r} \big)
    S_\omega^{d^{(q-1)\omega} - 1} \right)_r \\
  & = \left( R_\omega S_\omega^{d^{(q-1)\omega} - 1} \right)_{r} +
  (x^{[1]}_r - y^{[1]}_{r-1}), 
\end{split}
\end{equation*}
since $y_0^{[1]} = x_1^{[1]} = 1$.

Using equation \eqref{sum thing}, we
see that
\begin{equation*}
\begin{split}
  & x^{[q]}_r - y^{[q]}_{r-1} = \sum_{i=1}^r x^{[q-1]}_i \left(
    (R_\omega)^i (S_\omega)^{d^{q
        \omega} - i} \right)_{r} \\
  & - \sum_{j=0}^{r-1} x^{[q-1]}_{j+1} \left( (R_\omega)^{j+1}
    (S_\omega)^{d^{(q-1)\omega} - j - 1} \right)_{r} + (x^{[1]}_r -
  y^{[1]}_{r-1}) (x_1^{[q-1]}) \\
  & + (x^{[q-1]}_r - y^{[q-1]}_{r-1}) \left( (R_\omega)^r
    (S_\omega)^{d^{(q-1) \omega} - r} \right)_{r}  \\
  & = (x^{[1]}_r - y^{[1]}_{r-1}) (x_1^{[q-1]}) + (x^{[q-1]}_r -
  y^{[q-1]}_{r-1}),
\end{split}
\end{equation*}

We have $y^{[q-1]}_0 = x^{[q-1]}_1 = 1$, since $y^{[1]}_0 = x^{[1]}_1
= 1$.  Thus, assuming inductively that
$$ x^{[q-1]}_r - y^{[q-1]}_{r-1} =  (q-1)(x^{[1]}_r - y^{[1]}_{r-1}),$$
we have
\begin{equation}
x^{[q]}_r - y^{[q]}_{r-1} = q(x^{[1]}_r - y^{[1]}_{r-1} ).  
\end{equation}
Note in particular that $n_{q \omega} = r$ for all $q$, so $n_k$ is
bounded for all $k$, as desired.
 
Now,
$$
\lim_{q \to \infty} \frac{\log | q (x^{[1]}_r - y^{[1]}_{r-1}
  )|_v}{d^{q \omega}} = 0$$
and $\tau_{q \omega} = x^{[q]}_r -
y^{[q]}_{r-1}$.  Since $S_{q \omega} (1,0)$ is simply $y^{[q-1]}_0 =
1$, we have
$$\lim_{q \to \infty} \frac{\log|\tau_{q \omega}|_v}{d^{q \omega}} = 0
= \lim_{q \to \infty} \frac{\log |S_{q \omega}(0,1)|_v}{d^{q
    \omega}},$$
which give us \eqref{as desired} and thus completes
our proof.

\end{proof}

\subsection{Proofs of the main theorems}\label{proofs}
Now, we can show that the integral $\int_{\mathbb{P}^{1}(\bC_v)} \log
|t - \beta|_v d \mu_{\varphi,v}$ can be computed by taking the limit
of the average of $\log |\beta-w|_v$ on the points in
$\varphi^{-k}(\alpha)$, as $k \to \infty$, for any nonexceptional
point $\alpha$.

\begin{theorem}\label{backwards}
  Let $\alpha = [s:u]$ be a nonexceptional point in
  $\bP^1(\overline{K})$.  Then for any nonzero polynomial $F(t) \in {\overline
  K}[t]$ we have
\begin{equation*}
 \int_{\mathbb{P}^{1}(\bC_v)}
\log|F|_v \, d\mu_{\varphi,v}
  = \lim_{k \to \infty} \frac{1}{d^{k}}
  \sum_{\substack{\varphi^{k}([w:1]) = \alpha \\ F(w) \not= 0}}  \log | F(w) |_v.
\end{equation*}
where the $[w:1]$ for which
$\varphi^{k}([w:1]) = \alpha$ are counted with multiplicity.
\end{theorem}

\begin{proof}
  
  The polynomial $F$ factors as $F(t) = \gamma \prod_{i=1}^n
  (t-\beta_i)$ where $\gamma$ and $\beta_1, \dots, \beta_n$ are
  elements of ${\overline K}$.  For each $\beta_i$, the multiplicity
  of $\beta_i$ in $(\varphi^k)^*\alpha$ is at most $r((\varphi^k)^*
  \alpha)$ (where $r((\varphi^k)^* \alpha)$ is defined as in Section
  \ref{dio}).  Since $\alpha$ is not exceptional, we have $ \lim_{k
    \to \infty} \frac{r((\varphi^k)^* \alpha)}{d^k} = 0,$ by Lemma
  \ref{r(D)}.  Thus,
$$
\lim_{k \to \infty} \frac{1}{d^{k}}
\sum_{\substack{\varphi^{k}([w:1]) = \alpha \\ w \not= \beta_j}} \log |w -
\beta_j|_v = \lim_{k \to \infty} \frac{1}{d^{k}}
\sum_{\substack{\varphi^{k}([w:1]) = \alpha \\ F(w) \not= 0}} \log |w -
\beta_j|_v $$
for each $\beta_j$.  Hence, it suffices to show that
\begin{equation}\label{back eq}
\int_{\mathbb{P}^{1}(\bC_v)}
\log|t-\beta|_v \, d\mu_{\varphi,v} = 
  \lim_{{k} \to \infty} \frac{1}{d^{k}}
  \sum_{\substack{\varphi^{k}([w:1]) = \alpha \\ w \not= \beta}} \log | w - \beta |_v 
\end{equation}  
for any $\beta \in \bK$.

  Note that $\varphi^{k}([w:1]) = [s:u]$ if and only if $u P_{k}(w,1)
  - s Q_{k}(w,1) = 0$.  Thus, as polynomials in $t$, we have
$$   u P_{k}(t,1) - s Q_{k}(t,1) = \eta_{k} \prod_{\varphi^{k}([w:1]) =
  [s:u]} (t - w),$$
where $\eta_k \in \bK$.  We write
$$ u P_{k}(t,1) - s Q_{k}(t,1) = (t-\beta)^{w_k} G_k(t,1)$$
for a polynomial $G_k$ such that $G_k(\beta,1) \not= 0$, as in Proposition
\ref{bad back}.  Note that
$$ G_k(t,1) = \eta_{k} \prod_{\substack{\varphi^{k}([w:1]) =
  \alpha \\ w \not= \beta}} (t - w).$$

Plugging $\beta$ in for $t$ and taking logs of absolute
values gives
\begin{equation}
\log |G_k(\beta,1)|_v = \log |\eta_{k}|_v
+ \sum_{\substack{\varphi^{k}([w:1]) =
  [s:u] \\ w \not= \beta}} \log |w - \beta|_v.
\end{equation}

Applying Proposition \ref{bad back} therefore yields
\begin{equation}\label{almost}
\begin{split}  
  \lim_{{k} \to \infty} \frac{1}{d^{k}} &
  \sum_{\substack{\varphi^{k}([w:1]) = \alpha \\ w \not= \beta}}  \log |
  w - \beta |_v + \frac{\log | \eta_{k}
    |_v}{d^{k}} \\
  & = \lim_{{k} \to \infty} \frac{\log \max (|P_{k}(\beta,1)|_v,
    |Q_{k}(\beta,1)|_v) }{d^{k}}.
\end{split}
\end{equation}

Now, writing 
$$
u P_{k}(T_0,T_1) - s Q_{k}(T_0,T_1) = T_1^{w_k} V_k(T_0,T_1)$$
for
some polynomial $V_k$ such that $V_k(1,0) \not= 0$, we see that
$\eta_k = V_k(1,0)$.  Applying Proposition \ref{bad back}, we obtain
$$
\lim_{k \to \infty} \frac{\log|\eta_k|_v}{d^k} = \lim_{k \to
  \infty} \frac{\log \max (|P_{k}(1,0)|_v, |Q_{k}(1,0)|_v) }{d^{k}}.$$
Substituting this equality into \eqref{almost} gives
\begin{equation}
\begin{split}  
  \lim_{k \to \infty} \frac{1}{d^{k}}
  \sum_{\substack{\varphi^{k}([w:1]) = \alpha \\ [w:1]
      \not = \beta}}  \log | w - \beta |_v  & = \lim_{{k} \to \infty}
  \frac{\log \max (|P_{k}(\beta,1)|_v, |Q_{k}(\beta,1)|_v) }{d^{k}} \\
 & - \lim_{{k} \to \infty}
  \frac{\log \max (|P_{k}(1,0)|_v, |Q_{k}(1,0)|_v) }{d^{k}}.
\end{split}
\end{equation} 
Using Proposition \ref{more gen}, we obtain \eqref{back eq}. 
\end{proof}

Now, we show that the same result holds when we average $\log |\beta-w|_v$
over periodic points rather than inverse images of a point.

\begin{theorem}\label{periodic}
For any any polynomial $F \in {\overline K}[t]$ we have 
\begin{equation*}
\int_{\mathbb{P}^{1}(\bC_v)}
\log|F|_v \, d\mu_{\varphi,v} = 
  \lim_{{k} \to \infty} \frac{1}{d^{k}} \sum_{\substack{\varphi^{k}([w:1]) =
    [w:1] \\ F(w) \not= 0}} \log | F(w) |_v, 
\end{equation*}  
where the $[w:1]$ for which $\varphi^{k}([w:1]) = w$ are counted with
multiplicity.
\end{theorem}

\begin{proof}
  As in the proof of Theorem \ref{backwards}, it will suffice to show that
\begin{equation}\label{per eq}
\int_{\mathbb{P}^{1}(\bC_v)}
\log|t-\beta|_v \, d\mu_{\varphi,v} = 
  \lim_{{k} \to \infty} \frac{1}{d^{k}} \sum_{\substack{\varphi^{k}([w:1]) =
    [w:1] \\ w \not= \beta}} \log | w - \beta |_v
\end{equation}   
for any $\beta \in \bK$ (this follows from the fact that the
multiplicity of each $\beta_i$ as a $k$-periodic point is bounded for
all $k$ by Proposition \ref{bad per}).

We have $\varphi^{k}([w:1]) = [w:1]$ if and only if $P_{k}(w,1) - w
Q_{k}(w,1) = 0$.  Thus,
$$P_{k}(t,1) - t Q_{k}(t,1) = \gamma_{k} \prod_{\varphi^{k}([w:1]) =
  [w:1]} (t - w),$$
for some $\gamma_k \in \bK$.  We write
$$
P_{k}(t,1) - t Q_{k}(t,1) = (t-\beta)^{n_k} H_k(t,1)$$
for a
polynomial $H_k$ such that $H_k(\beta,1) \not= 0$.  We have
$$
H_k(t,1) = \gamma_k \prod_{\substack{\varphi^{k}([w:1]) = [w:1] \\ 
    w \not= \beta}} (t-w).$$

Then, plugging $\beta$ in for $t$, taking logs of absolute values, and
applying Proposition \ref{bad per} gives
\begin{equation}\label{almost per}
\begin{split}  
  \lim_{{k} \to \infty} \frac{1}{d^{k}}
  \sum_{\substack{\varphi^{k}([w:1]) = [w:1] \\ w \not= \beta}} & \log | \beta
  - w |_v + \frac{\log | \gamma_{k}
    |_v}{d^{k}} \\
  & = \lim_{{k} \to \infty} \frac{\log \max (|P_{k}(\beta,1)|_v,
    |Q_{k}(\beta,1)|_v) }{d^{k}}.
\end{split}
\end{equation}
Writing
$$ T_1 P_{k}(T_0,T_1) - T_0 Q_{k}(T_0,T_1) = T_1^{n_k} W_k(T_0,T_1)$$
for a polynomial $W_k$ such that $W_k(1,0) \not= 0$, we see that
$\gamma_k = W_k(1,0)$.  By Proposition \ref{bad per}, we have
$$
\lim_{k \to \infty} \frac{\log|\gamma_k|_v}{d^k} = \lim_{k \to \infty}
\frac{\log \max
  (|P_{k}(1,0)|_v, |Q_{k}(1,0)|_v) }{d^{k}}.$$
Combining this equality with \eqref{almost per} and Proposition
\ref{more gen} yields \eqref{per eq}.  
 
\end{proof}

We are now ready to prove the results regarding the computation of the
canonical height $h_\varphi(\beta)$.  First, we'll need a lemma.  Note
that the lemma does not follow directly from the work of
Call an Goldstine (\cite{CG}), since they only prove that in a fixed
number field, the local canonical heights sum to the global canonical
height.  What is required here is slightly different.

\begin{lemma}\label{for height}
  Let $\beta = [a:b]$ in $\bP^1(\bK)$.  Let $[a_1:b_1], \dots,
  [a_n:b_n]$ be the conjugates of $[a:b]$ under the action of
  $\Gal(\bK/K)$.  Then
\begin{equation}\label{lim sum}
\begin{split}
   [K(\beta) : K] & (\deg K ) h_\varphi([a:b]) \\
   & =  \sum_{\text{places $v$ of $K$}}
  \lim_{k \to \infty} \sum_{i=1}^n \frac{\log \max (|P_{k}(a_i,b_i)|_v,
    |Q_{k}(a_i,b_i)|_v) }{d^{k}}.
\end{split}
\end{equation}
\end{lemma}
\begin{proof}
For all but finitely many $v$, we have $|a_i|_v = |b_i|_v = 1$.
Furthermore, for all but finitely many $v$, we have 
\begin{equation}\label{good}
\log \max(|P_k(s,t)|_v, |Q_k(s,t)|_v) = 0
\end{equation}
for all $k$ whenever $|s|_v = |t|_v = 1$.  This is true, for example,
at all nonarchimedean $v$ of good reduction for $\varphi$ in the sense
of \cite{STP}.  Indeed, when $v$ is a finite place, \eqref{good} will
hold for all $|s|_v = |t|_v = 1$ unless either
$|\Res(P(T_0,1),Q(T_0,1))|_v$ or $|\Res(P(1,T_1),Q(1,T_1))|_v$ is less
than 1, where $\Res$ is the usual resultant of two polynomials (see
\cite[p.~279, Proposition 4]{BK}).  Thus, we can interchange the limit
and the sum on the right-hand side of \eqref{lim sum} so that
\begin{equation}\begin{split}\label{inter}
\lim_{k \to \infty} & \sum_{\text{places $v$ of $K$}}
   \sum_{i=1}^n \frac{\log \max (|P_{k}(a_i,b_i)|_v,
    |Q_{k}(a_i,b_i)|_v) }{d^{k}} \\
& = \sum_{\text{places $v$ of $K$}}
  \lim_{k \to \infty} \sum_{i=1}^n \frac{\log \max (|P_{k}(a_i,b_i)|_v,
    |Q_{k}(a_i,b_i)|_v)}{d^{k}}.
\end{split}
\end{equation}
Now, let $L$ be the field $K(\beta)$ and let $w$ be a place of $L$ that
extends the place $v$ of $K$; we write $w \mid v$.  The field $L$ has
$n$ embeddings $i: L \hookrightarrow \bC_v$; for exactly $[L_w:K_v]$
of these embeddings, we have $| i(x) |_v = |x|_w$ for all $x \in L$.
This yields $[L_w:K_v]$ conjugates $[a':b']$ of $[a:b]$ such that
$|a|_w = |a'|_v$ and $|b|_w = |b'|_w$.  Hence, we see that
\begin{equation*}
\begin{split}
  \sum_{i=1}^n & \log \max (|P_{k}(a_i,b_i)|_v, |Q_{k}(a_i,b_i)|_v) \\ 
  & = \sum_{w \mid v} [L_w:K_v] \log \max (|P_{k}(a,b)|_w,
  |Q_{k}(a,b)|_v).
  \end{split}
\end{equation*}  
Thus, we have
\begin{equation*}
\begin{split}
  \sum_{\text{places $v$ of $K$}} \sum_{i=1}^n \log \max
  (|P_{k}(a_i,b_i)|_v, |Q_{k}(a_i,b_i)|_v) \\
   = [K(\beta):K] (\deg K)
  h(\varphi^k([a:b])),
\end{split}
\end{equation*}
by \eqref{Weil}.  It follows from \eqref{can} and \eqref{inter} that
we therefore have
\begin{equation*}
\begin{split}
  \sum_{\text{places $v$ of $K$}} & \lim_{k \to \infty} \sum_{i=1}^n
  \frac{\log \max (|P_{k}(a_i,b_i)|_v,
    |Q_{k}(a_i,b_i)|_v)}{d^{k}}  \\
  & = [K(\beta):K](\deg K) \lim_{k \to \infty}
  \frac{h(\varphi^k([a:b]))}{d^k} \\
  & = [K(\beta):K](\deg K) h_\varphi([a:b]).
\end{split}
\end{equation*}
\end{proof}

\begin{theorem}\label{back height}
  Let $\alpha$ be any point in $\bP^1(\bK)$ that is not an exceptional
  point of $\varphi$.  Then, for any $\beta \in \bK$ and any nonzero
  irreducible $F \in K[t]$ such that $F(\beta) = 0$, we have
\begin{equation*}
\begin{split}
  (\deg K) & (\deg F) (h_\varphi(\beta)   - h_\varphi(\infty)) \\
  & =\sum_{\text{places $v$ of $K$}} \lim_{{k} \to \infty}
  \frac{1}{d^{k}} \sum_{\substack{\varphi^{k}([w:1]) = \alpha \\ F(w)
      \not= 0}} \log | F(w) |_v,
\end{split}
\end{equation*}
where the $[w:1]$ for which $\varphi^{k}([w:1]) = \alpha$
are counted with multiplicity.
\end{theorem}
\begin{proof}

Write $F(t) = \gamma \prod_{i=1}^n (t - \beta_i)$ where $\gamma \in K$
and the $\beta_i$ are the conjugates of $\beta$ under the action of
$\Gal(\bK/K)$.  By the product formula, we have $\sum_{\text{places
    $v$ of $K$}} \log |\gamma|_v = 0$.  Thus, using Theorem
\ref{backwards} and Proposition \ref{more gen}, we see that
\begin{equation}\label{above}
\begin{split}
  & \sum_{\text{places $v$ of $K$}} \lim_{k \to \infty}
  \frac{1}{d^{k}}
  \sum_{\substack{\varphi^{k}([w:1]) = \alpha \\ F(w) \not= 0}}  \log | F(w) |_v \\
  & = \sum_{\text{places $v$ of $K$}} \lim_{k \to \infty}
  \frac{1}{d^{k}}
  \sum_{\substack{\varphi^{k}([w:1]) = \alpha \\ F(w) \not= 0}}  \log \left| \prod_{i=1}^n (w - \beta_i)  \right|_v \\
  & = \sum_{i=1}^n \lim_{k
    \to \infty} \frac{\log \max (|P_{k}(\beta_i,1)|_v, |Q_{k}(\beta_i,1)|_v) }{d^{k}} \\
  & - (\deg F) \lim_{k \to \infty} \frac{\log \max (|P_{k}(1,0)|_v,
    |Q_{k}(1,0)|_v) }{d^{k}}.
\end{split}
\end{equation}
By Lemma \ref{for height}, the quantity on the last two lines is equal
to
$$ (\deg F)(\deg K)(h_\varphi(\beta) - h_\varphi(\infty)),$$
as desired.
\end{proof}

\begin{theorem}\label{per height}
For any $\beta \in \bK$ and any nonzero
  irreducible $F \in K[t]$ such that $F(\beta) = 0$, we have 
\begin{equation*}
\begin{split}
  (\deg K) & (\deg F) (h_\varphi(\beta)   - h_\varphi(\infty)) \\ 
  & =\sum_{\text{places $v$ of $K$}} \lim_{{k} \to \infty}
  \frac{1}{d^{k}} \sum_{\substack{\varphi^{k}([w:1]) = [w:1] \\ F(w)
      \not= 0}} \log | F(w) |_v,
\end{split}
\end{equation*}
where the $[w:1]$ for which $\varphi^{k}([w:1]) = w$ are counted with
multiplicity.
\end{theorem}
\begin{proof}
  The proof is the same as the proof of Theorem \ref{back height},
  using Theorem \ref{periodic} in place of Theorem
  \ref{backwards}.
\end{proof}




\section{A counterexample}\label{counterexample}
The main theorems of this paper are {\it not} true when we work over
the complex numbers $\bC$ rather than $\bK$.  Let $K = \bQ$ and let
$\varphi([x:y]) = [x^2:y^2]$ be the usual squaring map.  Let $v$ be
the archimedean place of $\bQ$, so that $\bC_v$ is just the usual
complex numbers $\bC$.  We define the function $\psi$ on the positive
integers recursively by $\psi(1) = 2$ and $\psi(n) = 2^{(n
  \psi(n-1))}$.  Let $\alpha = \sum_{n=1}^\infty 1/\psi(n)$ and let
$\beta = e^{2 \pi i \alpha}$.  Note that for any $t$, we have $|e^{2
  \pi i t} - 1| \leq \pi(t - [t])$, (where $[t]$ is the greatest
integer less than or equal to $t$).  Letting $\ell_n = \log_2
\psi(n)$, we then have
\begin{equation*}
\begin{split}
  & \frac{1}{2^{\ell_n}} \sum_{w^{2^{\ell_n}} = 1} \log |w - \beta|_v = \frac{\log |
  \beta^{\psi(n)}  - 1|}{\psi(n)} \\
  & \leq \frac{1}{\psi(n)} \log (\pi \left( \psi(n) \alpha - [\psi(n)
  \alpha]\right))
  \\
  & \leq \frac{1}{\psi(n)} \log \left(\pi \frac{\psi(n)}{\psi(n+1)}
  \sum_{j=0}^\infty \frac{1}{2^{j \psi(n+1)}}   \right) \\
  & \leq \log \pi + 1 - n \log 2 +  \log 2. 
\end{split}
\end{equation*}
Thus, $ \frac{1}{2^{\ell_n}} \sum_{w^{2^{\ell_n}} = 1} \log |\beta -
w|_v$ goes to $- \infty$ as $n \to \infty$, so
$$ \lim_{k \to \infty} \frac{1}{2^{k}} \sum_{w^{2^k} = 1} \log |w - \beta|_v$$
does not exist.  

\section{Applications and further questions}\label{applications}
\subsection{Lyapunov exponents}\label{lyap}
The Lyapunov exponent $L(\varphi)$ of a rational map
$\varphi:\bP^1_{\bC} \lra \bP^1_{\bC}$ (see \cite{mane2}) can be
defined as follows.  Choosing coordinates $[T_0:T_1]$ for
$\bP^1_{\bC}$, letting $t = T_0/T$, and writing $\varphi(t) =
P(t)/Q(t)$ for polynomials $P$ and $Q$, we define
$$
L(\varphi) = \int_{\bP^1(\bC)} \log | \varphi'(t) | d \mu_\varphi,$$
where $\mu_\varphi$ is the unique measure of maximal entropy
measure for $\varphi$ on $\bP^1$; this measure of maximal entropy is
the same as the Brolin-Lyubich measure discussed in
Section~\ref{brolin} (see \cite{Mane3}).  

The Lyapunov exponent can be computed via equidistribution on certain
subsequences of inverse images of nonexceptional points in
$\bP^1(\bC)$ (see \cite{demarco}, \cite{mane2}).  That is, given a
nonexceptional point $\alpha$ in $\bP^1(\bC)$, there is an infinite
strictly increasing sequence of integers $(m_i)_{i=1}^\infty$ such
that
$$
L(\varphi) = \lim_{i \to \infty} \frac{1}{(\deg \varphi)^{m_i}}
\sum_{\substack{\varphi^{m_i}(\beta) = \alpha \\ \varphi'(\beta) \not=
    0 \\ \beta \not= \infty}} \log |\varphi'(\beta)|.$$
It is not known,
however, if $L(\varphi)$ can be computed by taking the limit of the
average $\varphi'$ on the periodic points of $\varphi$.

When $\varphi$ is defined over a number field $K$, however, we obtain
the following result as a corollary of Theorem \ref{per eq}.  

\begin{cor}
  Let $K$ be a number field and let $\varphi:\bP^1_\bC \lra \bP^1_\bC$
  be a nonconstant rational map that is defined via base extension
  from a map $\varphi:\bP^1_K \lra \bP^1_K$.  Let $\varphi'$ be
  defined as above.  Then
$$ L(\varphi) = \lim_{k \to \infty} \frac{1}{(\deg \varphi)^k}
\sum_{\substack{\varphi^k(\xi) = \xi \\ \varphi'(\xi) \not= 0 \\ \xi
    \not= \infty}} \log |\varphi'(\xi)| . $$
\end{cor}
\begin{proof}
  We may write $\varphi'$ as a quotient of polynomials $A(t)/B(t)$
  with coefficients in $K$.  This yields $\log |\varphi'(t)| = \log
  |A(t)| - \log |B(t)|$.  The corollary then follows immediately from
  Theorem \ref{periodic}.
\end{proof}

This corollary says that if $\varphi$ is a rational function defined
over a number field, then the Lyapunov exponent of $\varphi$ is
completely determined by the derivative of $\varphi$ at the periodic
points of $\varphi$.  This means that the derivative of $\varphi$ at
the periodic points of $\varphi$ also determines the Hausdorff
dimension of the Julia set (see \cite{mane}).



 
\subsection{Symmetry of canonical heights}\label{symmetry pap}
In \cite{symmetry}, we show that when $\infty$ is not in the $v$-adic
Julia set of $\varphi$ for any archimedean $v$, we have
\begin{equation}\label{symmetry}
\lim_{{k} \to \infty} \frac{1}{d^{k}} \sum_{\varphi^{k}([w:1]) =
    [w:1]}  h(w) = \lim_{{\ell} \to \infty}
\frac{1}{2^{\ell}} \sum_{\xi^{{2}^\ell} = \xi }
h_\varphi(\xi). 
\end{equation}
This can be thought of as a symmetry relation, connecting $h$ of the
$\varphi$-periodic points with $h_\varphi$ of the roots of unity.  The
proof uses Theorem \ref{per height} along with Lyubich's
equidistribution theorem (\cite{Lyubich}) and some adelic intersection
theory (see \cite{zhangadelic} and \cite{zhangthese}).  We are also
able to use Theorem \ref{per height} to prove that
$$
h_\varphi(\beta) - h(\beta) \leq \lim_{{k} \to \infty}
\frac{1}{d^{k}} \sum_{\varphi^{k}([w:1]) = [w:1]} h(w) +
h_\varphi(\infty) + \log 2.$$

Our proof of \eqref{symmetry} does not work when $\infty$ is in the
$v$-adic Julia set of $\varphi$, for in that case the local height
${\hat h}_v$ is not bounded on the $v$-adic Julia set.  Unfortunately,
the $v$-adic Julia set is all of $\bP^1(\bC_v)$ when $v$ is
archimedean for many rational maps $\varphi$.  This is the case, for
example, when $\varphi$ is the map obtained by taking the
multiplication-by-2 map on an elliptic curve and modding out by the
hyperelliptic involution (such a map is called a Latt{\`e}s map).

On the other hand, the usual local height ${\hat h}_v(t)$ of an
element $t \in \bC_v$ is simply $\max (\log |t|_v, 0)$, which is only
a little bit different from $\log |t|_v$, and Theorem \ref{periodic}
proves a suitable equidistribution theorem for $\log |t|_v$.  We hope
to extend the techniques of this paper so that we can prove an analog
of Theorem \ref{periodic} for functions such as $\max (\log |t|_v,
0)$.

\subsection{Computing with points  of small height}\label{small}
The results in \cite{Bilu}, \cite{autissier}, \cite{BR}, \cite{FR1},
\cite{FR2}, and \cite{CL} all apply not only to the periodic points
and backwards iterates of a point that we treat in this paper but to
all points of small height in the algebraic closure of a number field
$K$.  For example, one the main theorems in \cite{BR}, \cite{FR1},
\cite{FR2}, and \cite{CL} states that for any continuous function $g$
on $\bP^1(\bC_v)$ and any infinite nonrepeating sequence of points $(\alpha_n)$
in $\bP^1(\overline{K})$ such that $\lim_{n \to \infty}
h_\varphi(\alpha_n) = 0$, one has
\begin{equation}\label{galois}
\lim_{n \to \infty} \frac{1}{|\Gal(\alpha_n)|}\sum\limits_{\sigma
  \in{\Gal(\alpha_n)}} g(\alpha_n^\sigma) = \int_ {\mathbb{P}^{1}(\bC_v)}
      g \, d\mu_{v,\varphi}, 
\end{equation}
where $\Gal(\alpha_n)$ is the Galois group of the Galois closure of
$K(\alpha_n)$ over $K$.

Baker, Ih, and Rumely (\cite{BIR}) and Autissier (\cite{aut-letter})
have produced counterexamples that show that \eqref{galois} does not
always hold when the function $g$ is replaced with $\log |F|_v$ for
$F$ a polynomial.  All of these examples involve infinite nonrepeating
sequences of points $(\alpha_n) \in \bQb$ such that $\lim_{n \to
  \infty} h(\alpha_n) = 0$ and
$$ \lim_{n \to \infty} \frac{1}{|\Gal(\alpha_n)|}\sum\limits_{\sigma
  \in{\Gal(\alpha_n)}} \log |\alpha_n^\sigma - 2| \not=
\int_{0}^{1} \log|e^{2\pi i \theta} - 2| d\theta.$$

The points $(\alpha_n)$ are not preperiodic in any of these examples
Thus, it may be possible to prove that the main results of this paper
continue to hold when we work with any nonrepeating sequence of Galois
orbits of preperiodic points.  This would imply the following
conjectured generalization of Siegel's theorem for integral points.

\begin{conjecture}[Ih]
  For any nonpreperiodic point $\beta \in \bP_{\fO_K}^1({\overline
    K})$, there are at most finitely many preperiodic points of
  $\varphi$ in $\bP_{\fO_K}^1(\overline{K})$ that are integral
  relative to $\beta$. (Here, $\fO_K$ is the ring of integers of $K$
  and $\alpha$ is said to be integral relative to $\beta$ if the
  Zariski closure of $\alpha$ does not meet the Zariski closure of
  $\beta$ in $\bP^1_{\fO_K}$.)
\end{conjecture}

Baker, Ih, and Rumely have proven that this is true when $\varphi$ is
a Latt{\`e}s map or the usual squaring map $x \mapsto x^2$ .  Using
Theorem \ref{per height} and arguing as in \cite{BIR} (or as in
\cite{SilSiegel}, which presents a related result), it is possible to
derive the following weak version of Ih's conjecture in general.

\begin{proposition}
  For any nonpreperiodic point $\beta \in \bP^1({\overline K})$, there
  are at most finitely many $n$ such that all $\alpha \in
  \bP^1({\overline K})$ of period $n$ are $\beta$-integral.
\end{proposition}


\providecommand{\bysame}{\leavevmode\hbox to3em{\hrulefill}\thinspace}
\providecommand{\MR}{\relax\ifhmode\unskip\space\fi MR }
\providecommand{\MRhref}[2]{%
  \href{http://www.ams.org/mathscinet-getitem?mr=#1}{#2}
}
\providecommand{\href}[2]{#2}

\end{document}